\documentclass[11pt,oneside,english]{amsart}
\usepackage[T1]{fontenc}
\usepackage[latin1]{inputenc}
\usepackage{geometry}
\usepackage{graphicx}
\usepackage{amssymb}
\usepackage{psfrag}

\makeatletter
\theoremstyle{plain}
\newtheorem{thm}{Theorem}[section]
\newtheorem{prop}[thm]{Proposition}
\newtheorem{lem}[thm]{Lemma}

\newtheorem{cor}[thm]{Corollary}
\theoremstyle{definition}

\newtheorem{prob}[thm]{Problem}
\newtheorem{exa}[thm]{Example}
\theoremstyle{remark}

\newcommand{\eqdef}{\stackrel{\rm def}{=}}
\newcommand{\xk}{X_1,\ldots,X_k}
\newcommand{\mk}{\mu^{(1)},\ldots,\mu^{(k)}}
\newcommand{\sk}{\sigma_1,\ldots, \sigma_k}
\newcommand{\lk}{\lambda^{(1)},\ldots,\lambda^{(k)}}
\usepackage{babel}

\newcommand{\Des}{\mathrm{Des}}
\newcommand{\maj}{\mathrm{maj}}
\newcommand{\Hilb}{\mathrm{Hilb}}
\newcommand{\St}{\mathcal {ST}}
\newcommand{\Asc}{\mathrm {Asc}}
\newcommand{\Codes}{\mathrm {Codes}}
\newcommand{\Coasc}{\mathrm {Coasc}}
\newcommand{\treuno}{\begin{picture}(12,6)
\put(1,0){\line(1,0){3}}
\put(1,3){\line(1,0){9}}
\put(1,6){\line(1,0){9}}
\put(1,0){\line(0,1){6}}
\put(4,0){\line(0,1){6}}
\put(7,3){\line(0,1){3}}
\put(10,3){\line(0,1){3}}
\end{picture}}
\newcommand{\duedue}{\begin{picture}(9,6)
\put(1,0){\line(1,0){6}}
\put(1,3){\line(1,0){6}}
\put(1,6){\line(1,0){6}}
\put(1,0){\line(0,1){6}}
\put(4,0){\line(0,1){6}}
\put(7,0){\line(0,1){6}}
\end{picture}}
\newcommand{\dueunouno}{\begin{picture}(9,6)
\put(1,-3){\line(1,0){3}}
\put(1,0){\line(1,0){3}}
\put(1,3){\line(1,0){6}}
\put(1,6){\line(1,0){6}}
\put(1,-3){\line(0,1){9}}
\put(4,-3){\line(0,1){9}}
\put(7,3){\line(0,1){3}}
\end{picture}}

\usepackage{amsmath,amssymb}
\usepackage{diagrams}
\makeatother
\begin{document}

\title[The refined multimahonian distribution]{Diagonal invariants and \\the refined multimahonian distribution}
\author{Fabrizio Caselli}

\maketitle

\begin{abstract}
Combinatorial aspects of multivariate diagonal invariants of the symmetric group are studied. As a consequence it is proved the existence of a multivariate extension of the classical Robinson-Schensted correspondence. Further byproduct are a pure combinatorial algorithm to describe the irreducible decomposition of the tensor product of two irreducible representations of the symmetric group, and new symmetry results on permutation enumeration with respect to descent sets.
\end{abstract}

\section{Introduction}
The invariant theory of finite subgroups generated by reflections has attracted  many mathematicians since their classification in the works of Chevalley \cite{C} and Shepard and Todd \cite{ST} with a particular attention on the combinatorial aspects of it. This is mainly due to the fact that the study of invariant and coinvariant algebras by means of generating functions leads naturally to nontrivial combinatorial properties of finite reflection groups. A crucial example in this context which is a link between the invariant theory and the combinatorics of the symmetric group is the Robinson-Schensted correspondence.
This correspondence (see \cite{R,S}) is a bijection between the symmetric group on $n$ elements and the set of ordered pairs of standard tableaux with $n$ boxes with the same shape. This is based on the row bumping algorithm and was originally introduced by Robinson to study the Littlewood-Richardson rule and by Schensted to study the lengths of increasing subsequences of a word.  This algorithm has found applications in the representation theory of the symmetric group, in the theory of symmetric functions and the theory of the plactic monoid.  Moreover, it is certainly fascinating from a combinatorial point of view and has inspired a considerable number of papers in the last decades. This correspondence has been generalized  to other Weyl groups, by defining ad hoc tableaux, or to semistandard tableaux in the so-called RSK-correspondence, by considering permutations as special matrices with nonnegative integer entries.

The main goal of this work is to exploit further the relationship between the Robinson-Schensted correspondence and the theory of invariants of the symmetric group.
By interpreting the Hilbert series with respect to a multipartition degree of certain (diagonal) invariant and coinvariant algebras in terms of (descents of) tableaux and permutations we deduce the existence of a multivariate extension of the Robinson-Schensted correspondence, which is based on the decomposition of tensor product of irreducible representations of the symmetric group. The idea of a relation between diagonal invariants and tensor product multiplicities for a finite subgroup of $GL(V)$ goes back to Solomon (see \cite[Remark 5.14]{So}) and pervades the results of Gessel \cite{Ge} on multipartite $P$-partitions. Although we can not define this correspondence explicitly, we can deduce from it an explicit combinatorial algorithm to describe the irreducible decomposition of the tensor product of two irreducible representations of the symmetric group. Finally, we show some further consequences in the theory of permutation enumeration.
\section{Background}
Let $V$ be a finite dimensional $\mathbb C$-vector space and $W$ be a finite subgroup of the general linear group $GL(V)$ generated by reflections, i.e. elements that fix a hyperplane pointwise. We refer to such a group simply as a \emph{reflection group}. The most significant example of such a group is the symmetric group acting by permuting a fixed linear basis of $V$. Other important examples are Weyl groups acting on the corresponding root space. In this paper we concentrate on the case of the symmetric groups (and some other related groups). Nevertheless, we preserve the symbol $W$ to denote the symmetric group $S_n$ on the $n$-element set $[n]\eqdef\{1,2,\ldots,n\}$.

Given a permutation $\sigma\in W$ we denote by
$$
\Des(\sigma) \eqdef \{i\in[n-1] :\, \sigma(i)>\sigma(i+1)\}
$$
the (right) \emph{descent set} of $\sigma$ and its \emph{major index} by
$$
\maj(\sigma) \eqdef \sum_{i\in \Des(\sigma)}i.
$$
For example if $\sigma=35241$ we have $\Des(\sigma)=\{2,4\}$ and $\maj(\sigma)=6$.
We recall the following equidistribution result  due to MacMahon (see \cite{M}).
\begin{thm}\label{macmahon}
 We have
\begin{eqnarray*}
 W(q) & \eqdef & \sum_{\sigma \in W}q^{\maj(\sigma)}=\sum_{\sigma \in W}q^{\mathrm{inv}(\sigma)}\\
 & = & \prod_{i=1}^n(1+q+q^2+\cdots+q^i),
\end{eqnarray*}
where $\mathrm{inv}(\sigma)=|\{(i,j):\,i<j \textrm{ and } \sigma(i)>\sigma(j)\}|$ is the number of inversions of $\sigma$.
\end{thm}

The dual action of a reflection group on $V^*$ can be extended to the symmetric algebra $S(V^*)$ of polynomial functions on $V$. If we fix a basis of $V$, the symmetric algebra is naturally identified with the algebra of polynomials $\mathbb C[X]$. Here and in what follows we use the symbol $X$ to denote an $n$-tuple of variables $X=(x_1,\ldots,x_n)$. The symmetric group $W$ acts on $\mathbb C[X]$ by permuting the variables. As customary, we denote by $\mathbb C[X]^W$ the ring of invariant polynomials (fixed points of the action of $W$). We also denote by $I^W_+$ the ideal of $\mathbb C[X]$ generated by homogeneous polynomials in $\mathbb C[X]^W$ of strictly positive degree. The \emph{coinvariant algebra} associated to $W$ is defined as the corresponding quotient algebra

$$R^W\eqdef\mathbb C[X]/ I^W_+.$$

The coinvariant algebra has important applications in the theory of representation since it is isomorphic to the group algebra of $W$ and in the topology of the flag variety since it is isomorphic to its cohomology ring.


If $R$ is a multigraded  $\mathbb C$-vector space we can record the dimensions of its homogeneous components via its Hilbert series

$$
\Hilb(R)(q_1,\ldots,q_k)\eqdef\sum_{a_1,\ldots,a_k\in\mathbb N} \dim (R_{a_1,\ldots,a_k})q_1^{a_1}\cdots q_k^{a_k},
$$
where $R_{a_1,\ldots,a_k}$ is the homogeneous subspace of $R$ of multidegree $(a_1,\ldots,a_k)$.\\
We note that, since the ideal $I^W_+$ is generated by homogeneous polynomials (by total degree) the coinvariant algebra is graded in $\mathbb N$. It turns out that the polynomial $W(q)$ appearing in Theorem \ref{macmahon} is the Hilbert series of the coinvariant algebra $R^W$:
\begin{equation}\label{wh}
 W(q)=\Hilb(R^W)(q).
\end{equation}

This is a crucial example of interplay between the invariant theory of $W$ and the combinatorics of $W$ (by Theorem \ref{macmahon}). All the other cases considered in this paper are algebraic and combinatorial variations and generalizations of this fundamental fact.

The coinvariant algebra affords also the structure of a multigraded vector space which refines the structure of graded algebra. This further decomposition  can be described in terms of descents of permutations and descents of tableaux and was originally obtained in a work of Adin, Brenti and Roichman \cite{ABR} for Weyl groups of type $A$ and $B$ (see also \cite{BC1} for Weyl groups of type $D$ and \cite{BB} for other complex reflection groups).

If $M$ is a monomial in $\mathbb C[X]$ we denote by $\lambda(M)$ its \emph{exponent partition}, i.e. the partition obtained by rearranging the exponents of $M$. We say that a polynomial is homogeneous of \emph{partition degree} $\lambda$ if it is a linear combination of monomials whose exponent partition is $\lambda$. We note that the exponent partition is not well-defined in the coinvariant algebra. For example, for $n=3$ the monomials $x_1^2$ and $x_2 x_3$ are in the same class in the coinvariant algebra (since $x^2_1-x_2x_3=x_1(x_1+x_2+x_3)-(x_1x_2+x_1x_3+x_2x_3)$), though they have distinct exponent partitions. Nevertheless the exponent partition will be fundamental in defining a ``partition degree'' also in the coinvariant algebra.

We recall the definition of the \emph{dominance order} in the set of partitions of $n$. We write $\mu\unlhd \lambda$, and we say that $\mu$ is smaller than or equal to $\lambda$ in the dominance order, if $\mu_1+\cdots +\mu_i\le \lambda_1+\cdots+\lambda_i$ for all $i$. We write $\mu\lhd \lambda$ if $\mu\unlhd \lambda$ and $\mu\neq \lambda$. We let $R^{(1)}_\lambda$ be the subspace of $R^W$ consisting of elements that can be represented as a linear combination of monomials with exponent partition smaller than or equal to $\lambda$ in dominance order. We also denote by $R^{(2)}_\lambda$ the subspace of $R^W$ consisting of elements that can be represented as a linear combination of monomials with exponent partition strictly smaller than $\lambda$ in dominance order. The subspaces $R^{(1)}_\lambda$ and $R^{(2)}_\lambda$ are also $W$-submodules of $R^W$ and we denote their quotient by
$$
R_\lambda \eqdef R^{(1)}_\lambda /R^{(2)}_\lambda.
$$
The $W$-modules $R_\lambda$ provide a further decomposition of the homogeneous components of the coinvariant algebra $R^W$.
\begin{thm}
There exists an isomorphism of $W$-modules
$$
\begin{diagram}
 \varphi:R^W_k   & \rTo^\cong  & \bigoplus_{|\lambda|=k }R_\lambda,
 \end{diagram}
$$
such that $\varphi^{-1}(R_\lambda)$ can be represented by homogeneous polynomials of partition degree $\lambda$.

\end{thm}

We can use this result to define a partition degree on the coinvariant algebra: we simply say that an element in $R_k^W$ is homogeneous of partition degree $\lambda$ if its image under the isomorphism $\varphi$ is in $R_\lambda$. 
We can therefore define the Hilbert polynomial of $R^W$ with respect to the partition degree by 
$$
\Hilb(R^W)(q_1,\ldots,q_n)=\sum_\lambda (\dim R_\lambda) q_1^{\lambda_{1}}\cdots q_n^{\lambda_{n}}.
$$

The dimensions of the $W$-modules $R_\lambda$ can be easily described in terms of descents of permutations. Given $\sigma \in W$ we define a partition $\lambda(\sigma)$ be letting
$$
\lambda(\sigma)_i=|\Des(\sigma)\cap \{i,\ldots,n\}|.$$
Note that the knowledge of $\lambda(\sigma)$ is equivalent to the knowledge of $\Des(\sigma)$ and that $\maj(\sigma)=|\lambda(\sigma)|$. Then we have the following result which can be viewed as a refinement of Equation \eqref{wh} (see \cite[Corollary 3.11]{ABR}).
$$\dim R_\lambda =|\{\sigma \in W: \lambda(\sigma)=\lambda\}|$$
Once we know the dimensions of the $W$-modules $R_\lambda$ we may wonder about their irreducible decomposition. For this reason we introduce  the \emph{refined fake degree polynomial} $f^\mu (q_1,\ldots, q_n)$  as the polynomial whose coefficient of $q_1^{\lambda_1}\cdots q_n^{\lambda_n}$ is the multiplicity of the representation $\mu$ in $R_\lambda$ if $\lambda$ is a partition, and zero otherwise, i.e. $$f^{\mu}(Q)=\sum_\lambda \langle \chi^{\mu},\chi^{R _\lambda} \rangle Q^\lambda,$$
where $Q^\lambda\eqdef q_1^{\lambda_1}\cdots q_n^{\lambda_n}$.
In this formula we denote by $\chi^{\rho}$ the character of a representation $\rho$ and by $\langle \cdot , \cdot \rangle$ the scalar product on the space of class functions on $W$ with respect to which the characters of the irreducible representations form an orthonormal basis. The polynomials $f^\mu (Q)$ have a very simple combinatorial interpretation based on standard tableaux that we are going to describe.\\
 Given a partition $\mu$ of $n$, the \emph{Ferrers diagram of shape $\mu$} is a collection of boxes, arranged in left-justified rows, with $\mu_i$ boxes in row $i$. A \emph{standard tableau} of shape $\mu$ is a filling of the Ferrers diagram of shape $\mu$ using the numbers from 1 to $n$, each occurring once, in such way that rows are increasing from left to right and columns are increasing from top to bottom. We denote by $\St$ the set of standard tableaux with $n$ boxes. For example the following picture

\vspace{3mm}
\begin{center}
\begin{picture}(80,80)
\put(0,80){\line(1,0){60}}
\put(0,60){\line(1,0){60}}
\put(0,40){\line(1,0){40}}
\put(0,20){\line(1,0){20}}
\put(0,0){\line(1,0){20}}
\put(0,0){\line(0,1){80}}
\put(20,0){\line(0,1){80}}
\put(40,40){\line(0,1){40}}
\put(60,60){\line(0,1){20}}
\put(7,7){7}
\put(7,27){4}
\put(7,47){2}
\put(7,67){1}
\put(27,47){6}
\put(27,67){3}
\put(47,67){5}
\put(-30,40){$T=$}
\end{picture}
\end{center}
\vspace{3mm}

\noindent  represents a standard tableau of shape $(3,2,1,1)$.
We say that $i$ is a \emph{descent} of a standard tableau $T$ if $i$ appears strictly above $i+1$ in $T$. We denote by $\Des(T)$ the set of descents of $T$ and we let $\maj(T)$ be the sum of its descents. Finally, we denote by $\mu(T)$ the shape of $T$. In the previous example we have $\Des(T)=\{1,3,5,6\}$ and so $\maj(T)=15$. As we did for permutations, given a tableau $T$ we define a partition $\lambda(T)$ by putting
$$(\lambda(T))_i=|\Des(T)\cap \{i,\ldots,n\}|.$$

It is known that irreducible representations of the symmetric group $W$ are indexed by partitions of $n$. We therefore use the same symbol $\mu$ to denote a partition or the corresponding Specht module. The following result appearing in \cite[Theorem 4.1]{ABR} describes explicitly the decomposition into irreducibles of the $W$-modules  $R_\lambda$ and refines a well-known result on the irreducible decomposition of the homogeneous components of $R^W$ attributed to Lusztig (unpublished) and to Kra\'skiewicz and Weyman (\cite{KW}).
\begin{thm}\label{fmu}The multiplicity of the representation $\mu$ in $R_\lambda$ is 
 $$
|\{T \textrm{ tableau}: \mu(T)=\mu \textrm{ and }\lambda(T)=\lambda\}|
$$ and so
$$
f^{\mu}(q_1,\ldots,q_n)=\sum_{\{T:\mu(T)=\mu\}}Q^{\lambda(T)}.
$$
\end{thm}
\section{Refined multimahonian distributions}
We let $\mathbb C[\xk]$ be the algebra of polynomials in the $nk$ variables $x_{i,j}$, with $i\in[k]$ and $j\in [n]$, i.e. we use the capital variable $X_i$ for the $n$-tuple of variables $x_{i,1},\ldots,x_{i,n}$.
We consider the natural action of  $W^k$  and of its diagonal subgroup $\Delta W$ on $\mathbb C[\xk]$. 
By means of the above decomposition of the coinvariant algebra we can also decompose the algebra
 $$\frac{\mathbb C[\xk]}{I^{W^k}_+}\cong \underbrace{R^W\otimes\cdots\otimes R^W}_k$$ 
in homogeneous components whose degrees are $k$-tuples of partitions with at most $n$ parts. In particular we say that an element in $\mathbb C[\xk]/I^{W^k}_+$ is homogeneous of \emph{multipartition degree} $(\lk)$ if it belongs to
$R^W_{\lambda^{(1)}}\otimes \cdots \otimes R^W_{\lambda^{(k)}}$ by means of the above mentioned canonical isomorphism. We are mainly interested in the subalgebra
$$
\left( \frac {\mathbb C[X_1,\ldots,X_k]}{I^{W^k}_+}\right)^{\Delta W}\cong \frac{\mathbb C[X_1,\ldots,X_k]^{\Delta W}}{J^{W^k}_+}.
$$
Here $J^{W^k}_+$ denotes the ideal generated by totally invariant polynomials with no constant term inside $\mathbb C[X_1,\ldots,X_k]^{\Delta W}$ and the isomorphism is due to the fact that the operator $$F\mapsto F^{\#}\eqdef \frac{1}{|W|}\sum_{\sigma \in \Delta W}\sigma(F)$$ is an inverse of the natural projection $ \mathbb C[X_1,\ldots,X_k]^{\Delta W}/J^{W^k}_+ \rightarrow  \left( \mathbb C[X_1,\ldots,X_k] / I^{W^k}_+\right)^{\Delta W} $. We can therefore consider the Hilbert polynomial 
$$
\Hilb\Big(\frac{\mathbb C[X_1,\ldots,X_k]^{\Delta W}}{J^{W^k}_+}\Big)\eqdef \sum_{\lk}\dim \Big(\frac{\mathbb C[X_1,\ldots,X_k]^{\Delta W}}{J^{W^k}_+}\Big)_{\lk} Q_1^{\lambda^{(1)}}\cdots Q_k^{\lambda^{(k)}}.
$$
In this formula the symbol $Q_i$ stands for the $n$-tuple of variables $q_{i,1},\ldots,q_{i,n}$ and the sum is over all partitions $\lk$ with at most $n$ parts.

Our next target is to describe the previous Hilbert series. For this we need to introduce one further ingredient.
We define the \emph{Kronecker coefficients} of $W$ by
\begin{eqnarray*}
 d_{\mk} & \eqdef & \frac{1}{|W|} \sum_{\sigma\in W} \chi^{\mu^{(1)}}(\sigma) \cdots \chi^{\mu^{(k)}}(\sigma)\\
 & = & \langle \chi^{\mu^{(1)}} \cdots \chi^{\mu^{(k-1)}},\chi^{\mu^{(k)}}\rangle_W,
\end{eqnarray*}
where $\mk$ are irreducible representations of $W$. In other words $d_{\mu^{(1)},\ldots,\mu^{(k)}}$ is the multiplicity of $\mu^{(k)}$ in the (reducible) representation $\mu^{(1)}\otimes\cdots \otimes \mu^{(k-1)}$.
These numbers have been deeply studied in the literature (see, i.e. \cite{BK,D,Re, Ro}) though they do not have an explicit description such as a combinatorial interpretation. A consequence of our main result is also a recursive combinatorial definition of the numbers $d_{\mu^{(1)},\ldots,\mu^{(k)}}$ which is independent of the character theory of $W$.\\
Now we can state the following result which relates the Hilbert series above with Kronecker coefficients and the refined fake degree polynomials.
\begin{thm}\label{fcgen}We have
\begin{eqnarray*}
\Hilb\Big(\mathbb C[\xk]^{\Delta W} / J^{W^k}_+\Big)(Q_1,\ldots,Q_k) & = & \sum_{\mk}d_{\mk}f^{\mu^{(1)}}(Q_1)\cdots f^{\mu^{(k)}}(Q_k)\\
 & = & \sum _{T_1,\ldots,T_k}d_{\mu(T_1),\ldots,\mu(T_k)} Q_1^{\lambda(T_1)}\cdots Q_k^{\lambda(T_k)}
\end{eqnarray*}
\end{thm}
\begin{proof}
The first equality is essentially due to Solomon (see \cite[Theorem 5.11]{So}). In fact, although Solomon's result
concerns only a simply graded G-module or a tensor power of a simply graded G-module (considered as a $G\wr S_n$ module), it can be easily generalized to the present context of a generic multigraded $G$-module.

The second equality follows directly from Theorem \ref{fmu}.
\end{proof}

We recall that the algebra $\mathbb C[\xk]^{\Delta W}$, being Cohen-Macauley (see \cite[Proposition 3.1]{Sta1}), is a free algebra over its subalgebra $ \mathbb C[\xk]^{W^k}$. This implies directly that if we consider $\mathbb C[\xk]$ as an algebra graded in $\mathbb N^k$ in the natural way, then,
\begin{equation}\label{pol}
\Hilb \Big(\mathbb C[\xk]^{\Delta W}/I^{W^k}_+\Big)(q_1,\ldots,q_k) = \frac{\Hilb (\mathbb C[\xk]^{\Delta W})(q_1,\ldots,q_k)}{\Hilb (\mathbb C[\xk]^{W^k})(q_1,\ldots,q_k)}.
\end{equation}
Now, the algebra $\mathbb C[\xk]$ is also multigraded by $k$-tuples of partitions with at most $n$ parts: we just say that a monomial $M$ is homogeneous of multipartition degree $(\lambda^{(1)},\ldots,\lambda^{(k)})$ if its exponent partition with respect to the variables $x_{i,1},\ldots x_{i,n}$ is $\lambda^{(i)}$ for all $i$. We write in this case $\lambda^{(i)}(M)\eqdef \lambda^{(i)}$ for all $i$ and $\Lambda(M)\eqdef(\lk)$. The refinement of Equation \eqref{pol} using the Hilbert series with respect to multipartition degree is no longer implied by the Cohen-Macauleyness of $\mathbb C[\xk]^{\Delta W}$. For this we need to use the existence of the representations $R_\lambda$ in a more subtle way. 
Given $\sigma \in W$ we define a monomial
$$
a_\sigma=\prod x_{\sigma(i)}^{\lambda(\sigma)_i}.
$$
By definition we clearly have $\lambda(a_{\sigma})=\lambda(\sigma)$.
In \cite{GS} and \cite{ABR} it is proved that the set of monomials $\{a_{\sigma}:\,\sigma \in W\}$ is a basis for the coinvariant algebra $R^W$. The proof in \cite{ABR} is based on a straightening law.  For its description we need to introduce an ordering on the set of monomials of the same degree: for $m$ and $m'$ monomials of the same total degree in $\mathbb C[X]$ we let $m\prec m'$ if
\begin{enumerate}
 \item $\lambda(m)\lhd \lambda(m')$; or
\item $\lambda(m)=\lambda(m')$ and $\textrm{inv}(\pi(m))>\textrm{inv} (\pi(m'))$,
\end{enumerate}
where $\pi(m)$ is the permutation $\pi$  having a minimal number of inversions such that the exponent in $m$ of $x_{\pi(i)}$ is greater than or equal to the exponent in $m$ of $x_{\pi(i+1)}$ for all $i$. 
The straightening law is the following: let $m$ be a monomial in $\mathbb C[X]$. Then $\mu:=\lambda(m)-\lambda(\pi(m))$ is still a partition and
\begin{equation}
 \label{straight} m=m_\mu \cdot a_{\pi(m)}+\sum_{m'\prec m}c_{m,m'}m',
\end{equation}
where $c_{m,m'}\in \mathbb C$ and $m_\mu$ is the monomial symmetric function. The straightening algorithm stated in \cite{ABR} uses elementary symmetric functions instead of monomial symmetric functions, but one can easily check that the two statements are equivalent. 
The fact that the set $\{a_{\sigma}:\,\sigma \in W\}$ is a basis of $R^W$ implies directly that the set of monomials
$$a_{\sk}\eqdef a_{\sigma_1}(X_1)\cdots a_{\sigma_k}(X_k)$$
is a basis for the coinvariant algebra of $W^k$, i.e. the algebra $\mathbb C[\xk]/I^{W^k}_+$.  Now, the monomials $a_{\sk}$ form a basis for the algebra $\mathbb C[X_1,\ldots X_k]$ as a free module over the subring $\mathbb C[\xk]^{W^k}$ of $W^k$-invariants (being a basis of the coinvariant algebra $\mathbb C[\xk]/I^{W^k}_+$), i.e.
$$\mathbb C[\xk]=\bigoplus_{\sk\in W}\mathbb C[\xk]^{W^k}a_{\sk}.$$ 
The following result states a triangularity property of this basis. If $(\mk)$ and $(\lk)$ are two $k$-tuples of partitions we write $(\mk)\unlhd (\lk)$ if $\mu^{(i)}\unlhd \lambda^{(i)}$ for all $i$ and we denote by $\mathbb C[\xk]_{\unlhd(\lk)}$ the space of polynomials spanned by monomials with multipartition degree $\unlhd (\lk)$. We similarly define $\mathbb C[\xk]_{\lhd(\lk)}$.
\begin{lem}\label{trian}
Let $M\in \mathbb C[\xk]$ be a monomial and let $$M=\sum_{\sk \in W}f_{\sk}a_{\sk},$$ where $f_{\sk}\in \mathbb C[\xk]^{W^k}$. Then this sum is restricted to those $\sk$ such that $\Lambda(M)-\Lambda(a_{\sk})$ is a $k$-tuple of partitions and
$$ f_{\sk}\in \mathbb C[\xk]_{\unlhd \Lambda(M)-\Lambda(a_{\sk})}.$$
\end{lem}
\begin{proof}
Given two monomials $M=m_1(X_1)\cdots m_k(X_k)$ and $M'=m_1'(X_1)\cdots m_k'(X_k)$  we let $M\prec M'$ if  $m_i\prec m'_i$ for all $i$. We proceed by a double  induction on the total degree and on $\prec$ within the set of monomials of the same multidegree. If $M$ has total degree zero the result is trivial. Otherwise let $\Lambda(M)=(\lk)$. If $M$ is minimal with respect to the ordering $\prec$ then $\lambda^{(i)}$ is minimal with respect to the dominance order for all $i$. If there exists $i$ such that $|\lambda^{(i)}|\geq n$ then $M=(x_{i,1}\cdots x_{i,n})M'$ and the result follows by induction since the total degree of $M'$ is strictly smaller than the degree of $M$. If $|\lambda^{(i)}|< n$ for all $i$ then $\lambda^{(i)}=(1^{k_i})$ by the minimality condition. Then $ m_i=a_{\pi(m_i)}$ for all $i$  and the result follows.
In the general case we can apply the straightening law (\ref{straight}) to all the $m_i$'s getting 
$$ M=f_M \cdot a_{\pi(m_1),\ldots,\pi(m_k)}+\sum_{M'\prec M}c_{M,M'}M',$$
where $f_M$ is a homogeneous $W^k$-invariant polynomial of multipartition degree $\Lambda(M)-\Lambda(a_{\pi(m_1),\ldots,\pi(m_k)})$. Then the result follows by induction.
\end{proof}
Now recall the already mentioned sequence of isomorphisms of $W$-modules
\begin{eqnarray*}
\frac{\mathbb C [\xk]^{\Delta W}}{J_+^{W^k}}&\cong &\left(\frac{\mathbb C [\xk]}{I_+^{W^k}}\right)^{\Delta W}\cong (\underbrace{R^W\otimes \cdots \otimes R^W}_{k\textrm{ times}})^{\Delta W}\\& \cong &\bigoplus_{\lk}(R_{\lambda^{(1)}}\otimes\cdots \otimes R_{\lambda^{(k)}})^{\Delta W}.
\end{eqnarray*}
Consider a basis of $(R_{\lambda^{(1)}}\otimes\cdots \otimes R_{\lambda^{(k)}})^{\Delta W}$. Every element of such a basis can be represented by a homogeneous element in  $\mathbb C [\xk]$ of multipartition degree $(\lk)$ (by definition) which is invariant for the action of $\Delta W$. In fact, if a representative $F$ of a basis element is not $\Delta W$-invariant we can consider its symmetrization 
$F^{\#}$ since, clearly,
$F$ and $F^{\#}$ represents the same class in $(R_{\lambda^{(1)}}\otimes\cdots \otimes R_{\lambda^{(k)}})^{\Delta W}$.
We denote by $\mathcal B(\lk)$ this set of representatives, i.e. $\mathcal B(\lk)$ is a set of polynomials in $\mathbb C [\xk]^{\Delta W}$ of multipartition degree $(\lk)$ whose corresponding classes form a basis of $(R_{\lambda^{(1)}}\otimes\cdots \otimes R_{\lambda^{(k)}})^{\Delta W}$. We denote by $\mathcal B$ the (disjoint) union of all $\mathcal B(\lk)$. By the Cohen-Macauleyness of $\mathbb C[\xk]^{\Delta W}$ we can deduce that the set $\mathcal B$ is a basis for $\mathbb C[\xk]^{\Delta W}$ as a free $\mathbb C[\xk]^{W^k}$-module (see \cite[Proposition 3.1]{Sta1}), i.e.
$$
\mathbb C[\xk]^{\Delta W}=\bigoplus_{b\in \mathcal B}\mathbb C[\xk]^{W^k} \cdot b.
$$ 
The following result implies a crucial triangularity property of the basis $\mathcal B$. 
\begin{lem}\label{verytrian}
Let $F\in \mathbb C[\xk]^{\Delta W}$ be homogeneous of multipartition degree $\Lambda(F)$. Then the unique expression
$$
F=\sum_{b\in \mathcal B}f_b b,
$$
with $f_b\in \mathbb C[\xk]^{W^k}$ for all $b\in \mathcal B$, is such that the sum is restricted to those $b\in\mathcal B$ for which $\Lambda(F)-\Lambda(b)$ is a $k$-tuple of partitions and
$$
f_b\in \mathbb C[\xk]^{W^k}_{\unlhd \Lambda(F)-\Lambda(b)}.$$
\end{lem}
\begin{proof}Let $\prec$ be a total order on the set of $k$-tuples of partitions of length at most $n$ satysfying the following two conditions
\begin{itemize}
 \item If $|\mu^{(1)}|+\cdots+|\mu^{(k)}|<|\lambda^{(1)}|+\cdots+|\lambda^{(k)}|$ then $(\mk)\prec(\lk)$;
\item If $\mu^{(i)}\unlhd \lambda^{(i)}$ for all $i$, then $(\mk)\prec(\lk)$.
\end{itemize}
We proceed by induction on the multipartition degree of $F$ with respect to the total order $\prec$. If $F$ has degree zero then the result is trivial. Otherwise let $\Lambda(F)=(\lk)$ be the multipartite partition degree of $F$. Then $F$ represents an element in $(R_{\lambda^{(1)}}\otimes\cdots \otimes R_{\lambda^{(k)}})^{\Delta W}$. Therefore,
$$F=\sum_{b\in \mathcal B(\lk)} c_b b$$
in  $(R_{\lambda^{(1)}}\otimes\cdots \otimes R_{\lambda^{(k)}})^{\Delta W}$. This means that 
$$F=\sum_{b\in \mathcal B(\lk)}c_b b +G$$
in $(\underbrace{R^W\otimes \cdots \otimes R^W}_{k\textrm{ times}})^{\Delta W}$, where $G$ is a $\Delta W$-invariant polynomial such that $$G\in \mathbb C[\xk]^{\Delta W}_{\lhd \Lambda(F)}.$$ Finally we deduce from this that
$$F=\sum_{b\in \mathcal B(\mk)}c_b b +G+H$$
in $\mathbb C[\xk]^{\Delta W}$, where $H$ belongs to $I_+^{W^k}$.
We can clearly assume that $G$ and $H$ are homogeneous with the same total multidegree of $F$. The induction hypothesis applies directly to $G$. Regarding $H$, by Lemma \ref{trian}, we can express
$H=\sum_{\sk}f_{\sk} a_{\sk}$ with $f_{\sk}\in \mathbb C[\xk]^{W^k}_{\unlhd \Lambda(F)-\Lambda (a_{\sk})}$ since $H$ is a sum of monomials of multipartition degree smaller than or equal to $\Lambda(F)$ in dominance order. Moreover, all the polynomials $f_{\sk}$ have positive degree since $H\in I_+^{W^k}$. Now we can apply the operator $\#$ to this identity and we get
$$H=\sum_{\sk}f_{\sk} a^{\#}_{\sk}.$$
Finally we can apply our induction hypothesis to the polynomials $a^{\#}_{\sk}$ since they have degree smaller than $F$ and the proof is completed by observing that, clearly,
$$
\mathbb C[\xk]_{\unlhd \Lambda}\cdot \mathbb C[\xk]_{\unlhd \Lambda'}\subseteq \mathbb C[\xk]_{\unlhd \Lambda+\Lambda'}.
$$
\end{proof}
We observe that Lemma \ref{verytrian} fails to be true for a generic homogeneous basis $\mathcal B$ of $\mathbb C[\xk]^{\Delta W}$ as a free $\mathbb C[\xk]^{W^k}$-module. We refer the reader to \cite[Section 4.2]{Rei}, \cite{A} and \cite{BL} for the explicit description of some bases of $\mathbb C[\xk]^{\Delta W}$ over $\mathbb C[\xk]^{W^k}$, with particular attention to the case $k=2$.

For notational convenience, if $\Lambda=(\lk)$ is a multipartition, we denote by 
$$\mathcal Q^\Lambda\eqdef Q_1^{\lambda^{(1)}}\cdots Q_k^{\lambda^{(k)}}=\prod_{i=1}^k\prod_{j=1}^n q_{i,j}^{\lambda^{(i)}_j}.$$
\begin{cor}\label{uou}
We have
$$
\Hilb\Big(\frac{\mathbb C[\xk]^{\Delta W}}{J^{ W^k}_+}\Big)(Q_1,\ldots,Q_k)= \frac{\Hilb(\mathbb C[\xk]^{\Delta W})(Q_1,\ldots,Q_k)}{\Hilb(\mathbb C[\xk]^{W^k})(Q_1,\ldots,Q_k)}=\sum_{b\in \mathcal B}\mathcal Q^{\Lambda(b)}$$ 
\end{cor}
\begin{proof}
The fact that $\Hilb\Big(\mathbb C[\xk]^{\Delta W} / J^{ W^k}_+\Big)(Q_1,\ldots,Q_k)=\sum_{b\in \mathcal B}\mathcal Q^{\Lambda(b)}$ is clear from the definition of the multipartition degree on $\mathbb C[\xk]^{\Delta W} / J^{ W^k}_+$ and the definition of the set $\mathcal B$. Lemma \ref{verytrian} implies that 
$$\dim \mathbb C[\xk]^{\Delta W}_{\unlhd \Lambda}= \sum_{b\in \mathcal B}\sum_{\{\Lambda':\Lambda'+\Lambda(b)\unlhd \Lambda\}}\dim \mathbb C[\xk]^{W^k}_{\Lambda'},$$
and similarly with $\lhd$ instead of $\unlhd$. Therefore
\begin{eqnarray*}
\dim \mathbb C[\xk]^{\Delta W}_{\Lambda} &=& \dim \mathbb C[\xk]^{\Delta W}_{\unlhd \Lambda}-\dim \mathbb C[\xk]^{\Delta W}_{\lhd\Lambda}\\
& = & \sum_{b\in \mathcal B}\sum_{\{\Lambda':\Lambda'+\Lambda(b)=\Lambda\}}\dim \mathbb C[\xk]^{W^k}_{\Lambda'}.
\end{eqnarray*}
Note that in the last sum there is only one summand corresponding to $\Lambda'=\Lambda-\Lambda(b)$ if this is a multipartition, and there are no summands otherwise. So we have
\begin{eqnarray*}
\Hilb(\mathbb C[\xk]^{\Delta W})& = & \sum_{\Lambda}\dim \mathbb C[\xk]^{\Delta W}_{\Lambda} \mathcal Q^\Lambda\\
& = & \sum_\Lambda \sum_{b\in \mathcal B}\sum_{\{\Lambda':\Lambda'+\Lambda(b)=\Lambda\}}\dim \mathbb C[\xk]^{W^k}_{\Lambda'}\mathcal Q^{\Lambda}\\
&=&\sum_{b\in \mathcal B}\sum_{\Lambda'}\dim \mathbb C[\xk]^{W^k}_{\Lambda'} \mathcal Q^{\Lambda'+\Lambda(b)}\\
&=&\sum_{b\in \mathcal B}\mathcal Q^{\Lambda(b)}\sum_{\Lambda'}\dim \mathbb C[\xk]^{W^k}_{\Lambda'}Q^{\Lambda'}\\
&=&\sum_{b\in \mathcal B}\mathcal Q^{\Lambda(b)} \Hilb (\mathbb C[\xk]^{W^k}).
\end{eqnarray*}
\end{proof}

Now we need to study the two Hilbert series of the invariant algebras $\mathbb C[\xk]^{\Delta W}$ and $\mathbb C[\xk]^{W^k}$ with respect to the multipartition degree. Before stating our next result we need to recall a classical theorem that can be attributed to Gordon \cite{G} and Garsia and Gessel \cite{GG} on multipartite partitions. 
We say that a collection $(f^{(1)},\ldots,f^{(k)})$ of $k$ elements of $\mathbb N^n$ is a $k$-partite partition if $f^{(i)}_j\geq f^{(i)}_{j+1}$ whenever $f^{(h)}_j=f^{(h)}_{j+1}$ for all $h<i$. For notational convenience we denote by $W^{(k)}\eqdef \{(\sk)\in W^k:\,\sigma_1\cdots\sigma_k=1)\}$. The main property of $k$-partite partitions that we need  is the following.
\begin{thm}\label{parpar}
There exists a bijection between the set of $k$-partite partitions and the set of $2k$-tuples $(\sigma_1,\ldots,\sigma_k, \mu^{(1)},\ldots, \mu^{(k)})$ such that
\begin{itemize}
\item $(\sk)\in W^{(k)}$;
\item $\mu^{(i)}$ is a partition with at most $n$ parts;
\item $\mu^{(i)}_j>\mu^{(i)}_{j+1}$ whenever $j\in \Des(\sigma_i)$.
\end{itemize}
The bijection is such that $\mu^{(i)}$ is obtained by reordering the coefficients of $f^{(i)}$.
\end{thm}
We can now prove the following formula for the quotient of the  Hilbert polynomials with respect to the multipartition degree associated to the invariant algebras of $\Delta W$ and  $W^k$.
\begin{thm} \label{GGgen}We have
$$\frac{\Hilb(\mathbb C[\xk]^{\Delta W})(Q_1,\ldots,Q_k)}{\Hilb(\mathbb C[\xk]^{W^k})(Q_1,\ldots,Q_k)}=
\sum_{(\sk)\in W^{(k)}}Q_1^{\lambda(\sigma_1)}\cdots Q_k^{\lambda(\sigma_k)} $$
\end{thm}
\begin{proof}
We observe that the set of monomials $X_1^{f^{(1)}}\cdots X_k^{f^{(k)}}$ as $(f^{(1)},\ldots,f^{(k)})$ varies among all possible $k$-partite partitions is a set of representatives for the orbits of the action of $\Delta W$ in the set of monomials in $\mathbb C[\xk]$.
By means of Theorem \ref{parpar} we can deduce that
$$\Hilb(\mathbb C[\xk]^{\Delta W})(Q_1,\ldots,Q_k)=\sum_{\substack{\sk,\\ \mk }}Q_1^{\mu^{(1)}}\cdots Q_k^{\mu^{(k)}},
$$
where the indices in the previous sum are such that they satisfy the conditions stated in Theorem \ref{parpar}.
We now observe that we have an equivalence of conditions
$$\mu^{(i)}_j>\mu^{(i)}_{j+1} \textrm{ whenever } j\in \Des(\sigma_i) \Longleftrightarrow  \mu^{(i)}-\lambda(\sigma_i) \textrm{ is a partition}.$$ Therefore we can simplify the previous sum in the following way
$$
\sum_{\substack{\sk,\\ \mk }}Q_1^{\mu^{(1)}}\cdots Q_k^{\mu^{(k)}} = \sum_{(\sk)\in W^{(k)}} \sum_{\nu^{(1)},\ldots,\nu^{(k)}}  Q_1^{\nu^{(1)}+\lambda(\sigma_1)}\cdots Q_k^{\nu^{(k)}+\lambda(\sigma_k)},
$$
where the last sum is on all possible $k$-tuples of partitions $\nu^{(1)},\ldots,\nu^{(k)}$ of length at most $n$. The result follows since, clearly, 
$$
\Hilb(\mathbb C[\xk]^{W^k} = \sum_{\nu^{(1)},\ldots,\nu^{(k)}}Q_1^{\nu^{(1)}}\cdots Q_k^{\nu^{(k)}}.$$
\end{proof}

Putting all these results together we obtain the following sequence of equivalent interpretations for what we may call the \emph{refined multimahonian distribution}.
\begin{thm}
We have
\begin{eqnarray*}
W(Q_1,\ldots,Q_k)& \eqdef & \sum_{T_1,\ldots,T_k}d_{\mu(T_1),\ldots,\mu(T_k)}Q_1^{\lambda(T_1)},\ldots,Q_k^{\lambda(T_k)}\\
                 & = & \sum_{\mk} d_{\mu^{(1)},\ldots,\mu^{(k)}} f^{\mu^{(1)}}(Q_1)\cdots f^{\mu^{(k)}}(Q_k)\\
                 & = & \Hilb \Big(\mathbb C[\xk]^{\Delta W} / J^{ W}_+)\Big)(Q_1,\ldots,Q_k)\\
                 & = & \frac{\Hilb(\mathbb C[\xk]^{\Delta W})(Q_1,\ldots,Q_k)}{\Hilb(\mathbb C[\xk]^{W^k})(Q_1,\ldots,Q_k)}\\
                 & = & \sum_{\sigma_1\cdots \sigma_k=1}Q_1^{\lambda(\sigma_1)}\cdots Q_k^{\lambda(\sigma_k)}
\end{eqnarray*}
\end{thm}
\begin{proof}
The four identities are the contents of Theorem \ref{fcgen}, Corollary \ref{uou} and Theorem \ref{GGgen} respectively.
\end{proof}
The cardinality of the set of $k$-tuples of permutations in $W^{(k)}$ having fixed descent sets was already studied by Gessel in \cite{Ge} and the idea to use Kronecker products (of quasi symmetric functions) is already present in his work. Nevertheless, Gessel's approach by means of quasi symmetric functions and the present approach my means of diagonal invariants are certainly different, and we believe that the present results can not be deduced from Gessel's work in a standard way.

The reason why we call the distribution $W(Q_1,\ldots,Q_k)$ refined is that one can consider its coarse version $W(q_1,\ldots q_k)$ obtained by putting $q_{i,j}=q_i$ for all $i$ and $j$. In this case one obtains the so-called \emph{multimahonian distribution} $\sum_{(\sk)\in W^{(k)}}q_1^{\maj(\sigma_1)}\cdots q_k^{\maj(\sigma_k)}$ which has been extensively studied in the literature (see, e.g.,  \cite{AR,BRS,BC,FS,GG}).
\begin{cor}\label{conj}
There exists a map $\mathcal T:W^{(k)}\longrightarrow \St^k$ satisfying the following two conditions:
\begin{enumerate}
\item For every $k$-tuple of tableaux $(T_1,\ldots,T_k)$, $$|\mathcal T^{-1}(T_1,\ldots,T_k)|=d_{\mu(T_1),\ldots,\mu(T_k)}.$$ In particular it depends only on the shapes of the tableaux $T_1,\ldots,T_k$;
\item  if $\mathcal T(\sigma_1,\ldots,\sigma_k)=(T_1,\ldots,T_k)$ then $\Des(T_i)=\Des(\sigma_i)$ for all $i=1,\ldots,k$.
\end{enumerate}
\end{cor}
The classical Robinson-Schensted correspondence (see \cite[\S 7.11]{Sta} for a description of this correspondence) provides a bijective proof of this corollary in the case $k=2$.

We can also conjecture that the correspondence $\mathcal T$ of Corollary \ref{conj} should be well-behaved with respect to cyclic permutations of the arguments in the sense that 
$$ \mathcal T(\sigma_1,\ldots, \sigma_k)=(T_1,\ldots,T_k) \Longrightarrow \mathcal T(\sigma_2,\ldots, \sigma_k, \sigma_1)=(T_2,\ldots,T_k, T_1).$$
\begin{prob}\label{rsmulti}
Find the map $\mathcal T$ of Corollary \ref{conj} explicitly.
\end{prob}
We observe that the resolution of this Problem would provide also an explicit combinatorial interpretation for the coefficients $d_{\mk}$.

\section{Combinatorial applications}
In this final section we deduce  some combinatorial results on Kronecker coefficients and permutation enumeration that follow from the results of the previous section. Our next goal is to show that we do not need to know the coefficients $d_{\mk}$ to solve Problem \ref{rsmulti}. This is because Corollary \ref{conj} uniquely determines the coefficients $d_{\mk}$ in the following sense.
\begin{prop}\label{duni}
Let  $\mathcal T:W^{(k)}\rightarrow \St^k$ be such that
\begin{enumerate}
 \item $|\mathcal T^{-1}(T_1,\ldots,T_k)|$ depends uniquely on the shapes of $T_1,\ldots,T_k$;
\item if $\mathcal T(\sigma_1,\ldots,\sigma_k)=(T_1,\ldots,T_k)$ then $\Des(T_i)=\Des(\sigma_i)$ for all $i=1,\ldots,k$.
\end{enumerate}
Then $\mathcal T$ satisfies the conditions of Corollary \ref{conj}, i.e. $|\mathcal T^{-1}(T_1,\ldots,T_k)|=d_{\mu(T_1),\ldots,\mu(T_k)}$ for all $(T_1,\ldots,T_k) \in \St^k$.
\end{prop}
This proposition is an immediate consequence of the following results which provides also an explicit entirely combinatorial algorithm to compute the coefficients $d_{\mk}$.

\begin{lem}
Let $\mu=(\mu_1,\ldots,\mu_r)$ be a partition of $n$ (with $\mu_r>0$). Then
\begin{enumerate}
 \item there exists a unique standard tableau $T_\mu$ of shape $\mu$ and descent set $\Des(T_\mu)=\{\mu_1,\mu_1+\mu_2,\ldots, \mu_1+\cdots+\mu_{r-1}\}$;
\item if $T$ is a standard tableau and $\Des(T)=\Des(T_\mu)$, then $\mu(T)\unrhd\mu$.
\end{enumerate}
\label{uiop}
\end{lem}
\begin{proof} The unique tableau $T_\mu$ satisfying these conditions is obtained by inserting the numbers form $1$ to $\mu_1$ in the first row of the Ferrers diagram of $\mu$, then the numbers from $\mu_1+1$ to $\mu_1+\mu_2$ in the second row and so on in the following rows. The uniqueness of $T_\mu$ is clear since the position of any single entry is forced by the descents conditions and the shape of the tableau.

 We observe that the number of descents of a tableau is at least the number of rows minus 1. This is because if $i\neq 1$ appears in the first column of $T$ then $i-1$ is necessarily above it and so $i-1\in \Des(T)$. Suppose $T$ is a standard tableau with the same descents as $T_\mu$. Consider the subtableau $T_i$ of $T$ consisting of the boxes filled with entries bounded by $\mu_1+ \cdots +\mu_i$. Then $T_i$ has exactly $i-1$ descents. So, if $j$ is the number of rows of $T_i$, we have $i-1\geq j-1$, i.e. $j\leq i$, by the previous observation. This implies that the numbers from 1 to $\mu_1+ \cdots +\mu_i$  appear all in the first $i$ rows of $T$ and the claim follows.
\end{proof}
The following corollary is a recursion satisfied by the coefficients $d_{\mk}$ and proves Proposition \ref{duni}.
\begin{cor}\label{algo}
Let $\mk$ be a multipartition and let $D_i=\Des(T_{\mu^{(i)}})$. Then
$$
d_{\mk}=\sum_{\substack {(\sk)\in W^{(k)}: \\ \Des(\sigma_i)=D_i}} 1- \sum_{\left\{\substack{T_1,\ldots,T_k\in \St:\, \Des(T_i)=D_i \textrm{ and}\\ \left(\mu(T_1),\ldots, \mu(T_k)\right) \rhd(\mk) }\right\}}d_{\mu(T_1),\ldots,\mu(T_k)}.
$$
\end{cor}
\begin{proof}
By Corollary \ref{conj} we have
$$
\sum_{\substack{T_1,\ldots,T_k\in \St :\\ \Des(T_i)=D_i}} d_{\mu(T_1),\ldots,\mu(T_k)}=\sum_{\substack {(\sk)\in W^{(k)}: \\ \Des(\sigma_i)=D_i}} 1. $$
Now the claim follows by Lemma \ref{uiop}.
\end{proof}

\begin{exa} Let $n=4$ and $k=3$. We  compute $d_{\treuno, \duedue, \dueunouno}$ by means of Corollary \ref{algo}. In this case we have $D_1=\{3\}, D_2=\{2\}, D_3=\{2,3\}$. We have to determine all tableaux having these descents sets. Now we observe that the unique tableau having descent set $D_1$,  is $T_{\treuno}$ and the unique tableau having descent set $D_3$ is $T_{\dueunouno}$. On the other hand there are two tableaux having descent set $D_2$ and these are $T_{\duedue}$ and $\begin{picture}(24,12)
\put(0,-4){\line(1,0){8}}
\put(0,4){\line(1,0){24}}
\put(0,12){\line(1,0){24}}
\put(0,-4){\line(0,1){16}}
\put(8,-4){\line(0,1){16}}
\put(16,4){\line(0,1){8}}
\put(24,4){\line(0,1){8}}
\put(2,-2){\tiny{$3$}}
\put(2,6){\tiny{$1$}}
\put(10,6){\tiny{$2$}}
\put(18,6){\tiny{$4$}}
\end{picture}$. So we have

$$
d_{\treuno, \duedue, \dueunouno}=\sum_{\substack {(\sigma_1,\sigma_2,\sigma_3)\in W^{(3)}: \\ \Des(\sigma_i)=D_i}} 1-d_{\treuno,\treuno,\dueunouno}
$$ 
The sum in the previous formula is $2$ since the only two triplets in the index set are $(1243,1423,1432)$ and $(2341,2413,2431)$. So we deduce that $
d_{\treuno, \duedue, \dueunouno}=2-d_{\treuno,\treuno,\dueunouno}$. Now we compute $d_{\treuno,\treuno,\dueunouno}$ by applying again Corollary \ref{algo}. In this case we have $D_1=D_2=\{3\}$ and $D_3=\{2,3\}$. By the previous observations on the tableaux having these descent sets we deduce that
$$
d_{\treuno,\treuno,\dueunouno}=\sum_{\substack {(\sigma_1,\sigma_2,\sigma_3)\in W^{(3)}: \\ \Des(\sigma_i)=D_i}} 1. 
$$
The last sum equals 1 since the only element in the index set is $(1342,1243,1432)$ and we can conclude that $d_{\treuno, \duedue, \dueunouno}=1$.
\end{exa}
To conclude this paper we prove some new results on permutation statistics that can be deduced from Corollary \ref{conj}. The first observation is a direct consequence of the symmetry of Kronecker coefficients with respect to their arguments. Let $D_1,\ldots,D_k\subseteq [n-1]$ and let $\pi$ be any permutation on $\{1,\ldots,k\}$. Then
\begin{equation}\label{sym}
|\{(\sigma_1,\ldots,\sigma_k)\in W^{(k)}:\,\Des(\sigma_i)=
D_i\}|=|\{(\sigma_1,\ldots,\sigma_k)\in W^{(k)}:\,\Des(\sigma_i)=D_{\pi(i)}\}|.
\end{equation}
This naturally leads to consider the following.
\begin{prob}Find a combinatorial bijective proof for the identity \eqref{sym}.
 \end{prob}
It is plausible that the resolution of this problem could be a first step towards the resolution of Problem \ref{rsmulti}.

The classical Robinson-Schensted correspondence allows one to prove some results on permutation enumeration and in particular on the bimohonian distributions (see for example \cite{DF,FS, BRS}). So, it is natural to ask how we can generalize these properties using the existence of the multivariate Robinson-Schensted correspondence.
If $X$ is either a permutation on $n$ elements or a tableau with $n$ entries we denote by 
\begin{eqnarray*}
\Codes(X)&\eqdef &\{i:n-i\in \Des(X)\};\\
\Asc(X) & \eqdef & [n-1]\setminus \Des(X);\\
\Coasc(X)& \eqdef & \{i:n-i\in \Asc(X)\}.
\end{eqnarray*}
The following is a multivariate generalization of a result of Foata and Sch\"utzenberger (\cite[Theorem 2]{FS})
\begin{prop}\label{descodes}
For all $I \subset \{1,...,k\}$ there exists an involution
$ F_I:W^{(k)} \rightarrow W^{(k)}$
such that, if $ F_I(\sk)=(\tau_1,...,\tau_k)$, then
$$
\Des(\sigma_i)=\left\{\begin{array}{ll} \Codes(\tau_i) & \textrm{if }i\in I;\\ \Des (\tau_i) & \textrm{otherwise.}\end{array}\right.
$$
\end{prop}
\begin{proof}
We recall (see \cite{FS, Schutzen}) that there exists a bijection $T\mapsto T^J$ on the set of standard tableaux of the same shape such that $\Des(T)=\Codes(T^J)$. Then the result follows immediately from Corollary \ref{rsmulti}.
\end{proof}
It is clear that we can substitute in Proposition \ref{descodes} $\Des$ with $\Asc$ and $\Codes$ with $\Coasc$ obtaining an analogous result. We can unify and generalize Equation \eqref{sym} and Proposition \ref{descodes} in the following statement.
\begin{thm}\label{dcac}
Fix $k$ subsets $D_1,\ldots,D_k \subseteq [n-1]$ arbitrarily. Then for any  integer sequence $0\leq i_1\leq i_2 \leq i_3 \leq k$ and for any permutation $\pi $ on $[k]$ the cardinality $C(i_1,i_2,i_3;\pi)$ of the set
$$
\left\{(\sk)\in W^{(k)}:\,\,D_i=\left\{ \begin{array}{ll}\Des(\sigma_{\pi(i)})& \textrm{if }0<i\leq i_1,\\ \Codes(\sigma_{\pi(i)})& \textrm{if }i_1 < i\leq i_2,\\  \Asc(\sigma_{\pi(i)})& \textrm{if }i_2< i \leq i_3, \\ \Coasc(\sigma_{\pi(i)})& \textrm{if }i_3<i\leq k, \end{array}
 \right.\right\}
$$ depends only on the parity of $i_2$, and in particular it does not depend on the indices $i_1, i_3$ and on the permutation $\pi$. \end{thm}
\begin{proof}
Since the knowledge of one of the three sets  $\Codes(\sigma)$, $\Asc(\sigma)$, $\Coasc(\sigma)$ is equivalent to the knowledge of the set $\Des(\sigma)$ it is clear that, by \eqref{sym}, the number $C(i_1,i_2,i_3;\pi)$ does not depend on $\pi$. So we can assume that $\pi=Id$. Then, from Proposition \ref{descodes}, we deduce that $C(i_1,i_2,i_3;Id)$ depends only on $i_2$. 

To prove that $C(i_1,i_2,i_3;\pi)$ depends only on the parity of $i_2$ we consider the following permutation on $W^{(k)}$
$$
(\sk)\mapsto (\sigma_1w_0,w_0\sigma_2, \ldots, \sigma_{2h-1}w_0,w_0\sigma_{2h}, \sigma_{2h+1},\ldots,\sigma_{k}),
$$
where $h$ is an integer such that $2h\leq k$ and $w_0=(n, n-1, \ldots, 1)$ is the top element of $W$. Observing that $\Des(\sigma)=\Coasc(\sigma w_0)=\Asc(w_0\sigma)$ for all $\sigma\in W$, this map proves bijectively that $C(k,k,k,Id)=C(k-2h,k-2h,k-h,\pi)$, where $\pi=(2h+1,2h+2,\ldots,k,2,4,\ldots,2h, 1,3,\ldots, 2h-1)$ and $C(k-1,k-1,k-1,Id)=C(k-1-2h,k-1-2h,k-1-h),\pi')$, where $\pi'=(2h+1,2h+2,\ldots,k-1,2,4,\ldots,2h, 1,3,\ldots, 2h-1,k)$) and the proof is complete.
\end{proof}
As an example on how we can obtain properties on the Kronecker coefficients by means of Proposition \ref{duni} avoiding their definition in terms of characters of the symmetric group we show the following result. If $\lambda$ is a partition we denote by $\lambda'$ the partition conjugate to $\lambda$.
\begin{cor}
We have
$$
d_{\lk}=d_{ {\lambda^{(1)'}},\lambda^{(2)'},\lambda^{(3)},\ldots,\lambda^{(k)}}
$$
\end{cor}
\begin{proof}
It is clear that transposition $T \mapsto T^t$ is a bijection between standard tableaux of conjugate shapes such that $\Des(T)=\Asc(T^t)$. Then, by Corollary \ref{rsmulti} and Theorem \ref{dcac} we have
\begin{eqnarray*}
 \sum_{\substack{(T_1,\ldots,T_k)\in \St^k:\\ \Des(T_i)=D_i}}d_{\lambda(T_1)',\lambda(T_2)',\lambda(T_3)\ldots\lambda(T_k)} & = & \sum_{\substack{(T_1,\ldots,T_k)\in \St^k:\\ \Asc(T_i)=D_i,\,i=1,2\\ \Des(T_i)=D_i,\,i\geq3}}d_{\lambda(T_1),\ldots,\lambda(T_k)}
  =  \sum_{\substack{(\sk)\in W^{(k)}:\\ \Asc(\sigma_i)=D_i,\,i=1,2\\ \Des(\sigma_i)=D_i,\,i\geq3}}1\\
 & = & \sum_{\substack{(\sk)\in W^{(k)}:\\ \Des(\sigma_i)=D_i}}1 =  \sum_{\substack{(T_1,\ldots,T_k)\in \St^k:\\ \Des(T_i)=D_i}}d_{\lambda(T_1),\ldots, \lambda(T_k)}.
\end{eqnarray*}
Now we can conclude recalling the uniqueness of the coefficients $d_{\lk}$ given by Proposition \ref{duni}.\\
We note that this result can also be easily obtained from the definition of the Kronecker coefficients recalling that the irreducible representation indexed by $\lambda'$ is isomorphic to the tensor product of the irreducible representation indexed by $\lambda$ with the alternant representation.
\end{proof}

We conclude this paper with the following result that we state without proof. The Robinson-Schensted correspondence can be naturally generalized to the so-called RSK-correspondence (where the ``K'' stands for Knuth) between multisets with support in $\mathbb N^2$ of cardinality $n$ and pairs of semistandard tableaux of the same shape. One can show that  an analogous correspondence exists also in the multivariate case between multisets with support in $\mathbb N^k $ of cardinality $n$ and $k$-tuples of semistandard Young tableaux of size $n$ preserving ``descents''. Here the descent sets $\Des_1(A),\ldots\Des_k(A)$ of a multiset $A$ are by definition the descent sets of the permutations appearing in Theorem \ref{parpar} for the $k$-partite partition obtained reordering the elements of $A$. The descents of a semistandard tableau are the descents of its standardized tableau. 
\begin{thm}
There exists a correspondence $\mathcal T$ between multisets with support in $\mathbb N^k $ of cardinality $n$ and $k$-tuples of semistandard Young tableaux of size $n$ such that,
\begin{enumerate}
\item $|\mathcal T^{-1}(T_1,\ldots,T_k)|=d_{\mu(T_1),\ldots,\mu(T_k)}$;
\item if $\mathcal T(A)=(T_1,\ldots,T_k)$ then the multiplicity of $i$ in $T_j$ is equal to the multiplicity of $i$ within the $j$-th coordinates of the elements of $A$ and
\item $\Des(T_i)=\Des_i(A)$.
\end{enumerate}

\end{thm}
{\bf Acknowkedgements.} A particular gratitude goes to Vic Reiner for his precious comments on a preliminary draft of this paper.

{\sc Dipartimento di matematica, Universit\`a di Bologna, \\Piazza di Porta San Donato 5, \\ Bologna 40126, Italy}\\
\emph{E-mail address: }{\tt caselli@dm.unibo.it}
\end{document}